\documentclass[10pt,reqno]{amsart}

\usepackage[cp1251]{inputenc}

\usepackage{amssymb
,epsfig
}
\usepackage{graphicx}

\newtheorem{thm}{Theorem}
\newtheorem{prp}{Proposition}
\newtheorem{lemma}{Lemma}

\newtheorem{ex}{Example}
\newtheorem{rem}{Remark}
\newtheorem{alg}{Algorithm}
\begin{document}

\title{Coulomb control of polygonal linkages}

\author{G.Khimshiashvili*, G.Panina$^\dag$, D.Siersma$^\ddag$}

\address{*Ilia State University, Tbilisi,  Georgia,
e-mail: giorgi.khimshiashvili@iliauni.edu.ge
 $^\dag$ Institute for Informatics and Automation, St. Petersburg, Russia,
Saint-Petersburg State University, St. Petersburg, Russia,
e-mail:gaiane-panina@rambler.ru. $^\ddag$ University of Utrecht,
Utrecht, The Netherlands, e-mail: D.Siersma@uu.nl. }


\begin{abstract}
Equilibria of polygonal linkage with respect to Coulomb potential of
point charges placed at the vertices of linkage are considered. It
is proved that any convex configuration of a  quadrilateral linkage
is the point of global minimum of Coulomb potential for appropriate
values of charges of vertices. Similar problems are treated  for the
equilateral  pentagonal linkage. Some corollaries and applications
in the spirit of control theory are also presented.
\end{abstract}

\maketitle

\section{Introduction}
\subsection*{Preliminary remarks }
Motivated by the famous Maxwell conjecture on equilibria of point
charges \cite{max} (cf. also \cite{ganosh}) we deal with the the
Coulomb potential of a system of point charges placed at the
vertices of a (flexible) planar polygonal linkage.  We consider
Coulomb potential of a vertex-charged linkage as a meromorphic
function on the planar moduli space of the linkage and investigate
its critical points.

The scenario we have in mind is suggested by some recent research
concerned with the control of nanosystems \cite{kundernac}. As an
abstract analog of a real physical situation we suggest the
following setting which can be described as {\it controlling the
shape of linkage by the values of charges at its vertices}. The
basic implicit assumption and motivation is that a vertex-charged
linkage subject only to Coulomb interaction of its charged vertices
should take the shape with the minimal Coulomb potential.

This setting suggests several aspects and problems. In the present
paper, we concentrate on the following scenario. Given a planar
configuration of linkage we wish to find the vertex charges such
that the global minimum of the arising Coulomb potential is achieved
at the given configuration. Such a collection of charges will be
said to {\it stabilize} the given configuration. If any
configuration of linkage has a stabilizing system of charges we will
say that this linkage admits a {\it complete Coulomb control}.

 Assuming
that a stabilizing collection of charges exists what does the set of
all stabilizing charges look like? How many minima and other
critical points has a system of stabilizing charges?

An interesting special case of the previous problem arises if one
asks if any convex configuration can be stabilized by a system of
charges of the same sign.
 This problem will be called {\it Coulomb
control of convex configurations}.

 It should be noted that this
research arose as a natural continuation of our previous joint
results on Morse functions on moduli spaces of polygonal linkages
\cite{khipan}, \cite{khsi}, \cite{kpsz}. The present paper was
completed during a "Research in Pairs" session in CIRM (Luminy) in
January of 2013. The authors acknowledge the hospitality and
excellent working conditions at CIRM.

\subsection*{Definitions and results }

A \textit{polygonal linkage} $L$ is  defined by a collection of
positive numbers $l = (l_1,...,l_n)$, called {\it sidelengths},
which we express by writing $L = L(l)$.

Physically, a polygonal linkage is a collection of rigid bars of
lengths $l_i$ joined in a cycle by revolving joints. It is a
flexible mechanism which can admit different shapes,  with or
without intersections.

By $M(L)$ we denote the \textit{moduli space} of planar
configurations,
 that is, the space of all polygons with the prescribed edge lengths factorized
 by isometries of $\mathbb{R}^2$:

$$M(L)=\{(p_1,...,p_n)| p_i \in \mathbb{R}^2, |p_ip_{i+1}|=l_i, \ |p_np_{1}|=l_n\}/Iso(\mathbb{R}^2).$$

This is not exactly the moduli space $\mathcal{M}(L)$ treated in
\cite{khipan} and \cite{F}, where the space of polygons is
factorized by orientation preserving isometries. However, there is a
two-fold covering $\mathcal{M}(L)\rightarrow M(L).$

 By $M^C(L)$ we denote the set of all convex  configurations.
 We allow here non-strictly convex polygons,  that is, those having (at least) one angle
 equal to
$\pi$. The latters obviously form the boundary $\partial M^C(L)$.
The set of all strictly convex configurations (all angles are less
than $\pi$) is the interior $IntM^C$. It is known (see
\cite{GaianeModuliComb}) that $M^C(L)$ is  homeomorphic to a ball.
In this paper we only deal with  $n=4, \ 5$.  For a 4-bar polygonal
linkage, $M(L)$ is a (topological) circle, whereas $M^C(L)$ is
homeomorphic to a segment. For a 5-bar polygonal linkage, $M(L)$ is
(generically) a surface, whereas $M^C(L)$ is the disk $D^2$.

\bigskip

Putting a collection of charges $q_i$ at the vertices $p_i$ of a
configuration and considering the Coulomb potential of these charges
we get a  function defined on $M(L)$. We will refer to this setting
by speaking of a {\it vertex-charged linkage} with the system of
charges $q_i$.

 Recall that
the Coulomb potential $\tilde{E}$ of a system of point charges $ q_i
\in \mathbb{R}$ placed at the points $p_i$ of Euclidean plane
$\mathbb{R}^2$ is defined as

\begin{equation} \label{coulomb}
\tilde{E} = \sum_{i\neq j} \frac{q_iq_j}{x_{ij}},
\end{equation}
where $x_{ij} = \vert p_i - p_j \vert$ is the distance between $i$th
and $j$th charges.

Since we are only interested in critical points of Coulomb
potential, addition of a constant makes no difference. By the very
definition of polygonal linkage,  the distances corresponding to two
consequent vertices in formula (\ref{coulomb}) remain the same for
all configurations of linkage. Hence their sum is constant for any
fixed collection of charges and for our purposes it is sufficient to
work with the {\it effective} Coulomb potential $E$ of configuration
defined as

\begin{equation} \label{linkcoulomb}
E = \sum \frac{q_iq_j}{x_{ij}},
\end{equation}

 where $x_{ij}$ is the length
of diagonal between (non-neighboring) $i$th and $j$th vertex of the
configuration. We say that a collection of charges
\textit{stabilizes} the configuration $P$ if $E$ attains at $P$ its
minimal value. In this case we say that $P$ is \textit{the minimum
point} of $E$.

We explicate now the setting and notation. For $n=4$, we put one
positive charge $t$  at the first vertex. The rest three
vertex-charges are $+1$.

For $n=5$, we put two positive charges $s$ and $t$ at  any two
non-neighboring vertices and say that $s$ and $t$ are
\textit{controlling charges}. The rest three charges, called {\it
non-controlled charges}, are again equal to $+1$.

\begin{rem}
If all the pairs of consecutive sidelengths are different, $E$ is a
smooth function without poles. If not, in our setting we have only
the poles with positive residues, which do not affect our study of
minima points of $E$.
\end{rem}

\newpage

Our \textbf{main results} are:

\begin{enumerate}
    \item For any 4-bar linkage, and any positive charge $t$, there exists a unique  convex
    configuration which is a critical point of $E$. This is the global minimum of $E$, and it depends continuously
    on $t$.
    \item For any 4-bar linkage, and each convex configuration $P$, there exists a unique  charge $t$
     which (together with the non-controlled charges) stabilize $P$. In this case, $t$ is positive, and $P$ is the global minimum point.
    \item We have a \textit{complete Coulomb control for convex quadrilaterals}.
    This means a two-step navigating algorithm  bringing any 4-bar
configuration to an in
    advance prescribed convex position  ruling by a positive charge $t$.
    \item For any  convex  equilateral pentagon $P$, there exists a unique pair of positive charges $(s,t)$
    for which $P$ is a critical point of $E$. However, it is unclear whether  $P$ is the global minimum point.
\end{enumerate}

\section{Coulomb control problem of convex quadrilaterals}

We begin with considering Coulomb control for convex configurations
of a non-degenerate $4$-bar linkage $L$ with one positive
controlling charge $t>0$ at the first vertex and three equal charges
$+1$ at the other three vertices. For a convex configuration  $P$
 of $L$, we denote by $x$ and $y$
the lengths of its two diagonals, and by $E$ its (effective) Coulomb
potential with controlling charge $t$:

\begin{equation} \label{coulombfour}
E =\frac{t}{x} +\frac{1}{y}.
\end{equation}

\begin{lemma}\label{UniquetForFourGon}
For a given convex quadrilateral  $P \in M(L)$, there exists a
unique $t$ such that $P$ is a critical point of $E$ on $M(L)$. In
this case, $t$ is positive.
\end{lemma}
Proof.  In a neighborhood of $P$ in $M(L)$, we have a relation of
the form $y = y(x)$. Then the condition that $P$ is a critical point
of $E$ is
$$y^2/x^2=-ty'(x),$$
which defines $t$ uniquely.  For a convex configuration $P$,  each
flex of $P$ which increases one of diagonals, shortens the other
one. This means that $y'(x)$ is negative. Hence $t$ is positive.\qed

\begin{prp}\label{thm4gon}
For a 4-bar linkage and a positive charge $t$,
\begin{enumerate}
    \item $E$ has a unique  minimum  in the interior of $M^C(L)$
    (that is, among strictly convex configurations),
which depends on $t$ continuously.

    \item  $E$ has no critical points among non-convex non-intersecting 4-gons.

    \item There is at least one critical point of $E$, which is a self-intersecting 4-gon.
\end{enumerate}
\end{prp}
Proof.

(1) First, we show that all the points $(x^2,y^2)$  lie on an
elliptic curve.

\begin{figure}[h]
\centering
\includegraphics[width=10 cm]{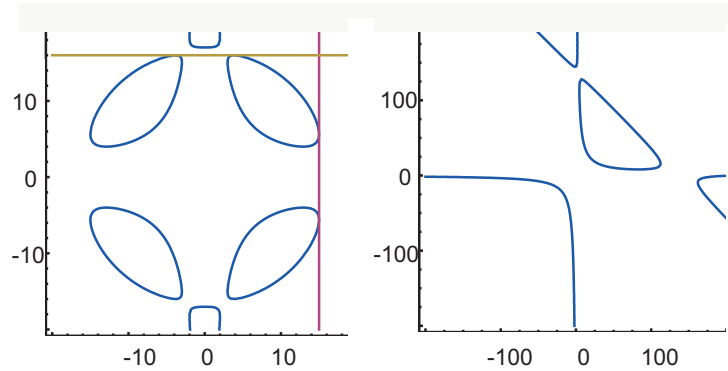}
\caption{ Diagonal relation in $x,y$ (left) and in $x^2,y^2$
(right)}\label{oval}
\end{figure}

Indeed, let  $$M(x,y)=\left(\begin{array}{ccccc}
0 & 1 & 1 & 1 & 1 \\
1 & 0 & a^2 & x^2 & d^2 \\
1 & a^2 & 0 & b^2 & y^2 \\
1 & x^2 & b^2 & 0 & c^2 \\
1 & d^2 & y^2 & c^2 & 0
\end{array}\right).$$
Put
$$ F(x,y) = \det(M(x,y))= -2 x^4 y^2-2 x^2y^4
 + 2x^2(d^2-a^2)(b^2-c^2 )+ 2y^2(b^2-a^2)(d^2-c^2) $$$$-2(ac-bd)(ac+bd)(a^2-b^2+c^2-d^2), $$
which is the  Cayley-Menger determinant. The equation  $ F(x,y) =0$
defines a real algebraic curve of degree 6, which depends only on
$x^2$ and $y^2$. In fact it is a cubic $g(w,z)$ with $w=x^2, \
z=y^2$.  The configuration space $M(L)$ maps bijectively to the
compact convex component (\textit{oval}) of the curve (inside a box
bounded by squared maximal length of diagonals), see Fig.
\ref{oval}. This component has has not more than two intersection
points with any straight line.

\begin{figure}[h]
\centering
\includegraphics[width=10 cm]{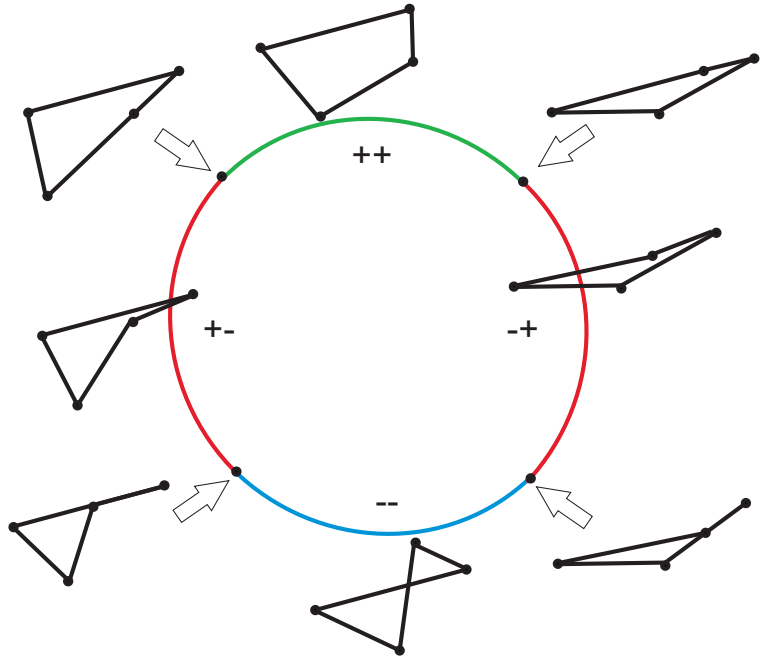}
\caption{ The configuration space  $M(L)$ is divided into four
parts}\label{Circle}
\end{figure}
Consequently, the oval contains exactly 2 points where
$\frac{\partial F}{\partial x} = 0$ and two other points where
$\frac{\partial F}{\partial y} = 0$. These points divide the oval
into 4 segments, each determined by the signs of partial derivatives
 (see Fig. \ref{Circle}).
\begin{enumerate}
    \item Convex 4-gons (that is, $M^C$)
correspond to the "$++$" case. A vertex of  configuration is
\textit{aligned}  if the angle adjacent to it is $\pi$. For a
non-degenerate 4-bar linkage, there exist exactly 2 vertices  which
can be aligned. The corresponding configurations are the boundary
points of $M^C$.
    \item Non-convex and non-self-intersecting 4-gons correspond to "$+-$" and "$-+$" parts.
    \item Self-intersecting 4-gons correspond to the "$- -$" part.
\end{enumerate}

We prove now that, for a convex $P$, we have $y''_{xx}<0$.

Indeed, locally $w$ is a function in $z$. We have  $w=w(z)$, therefore $y^2=w(x^2)$, which implies
$$y=w^{1/2}(x^2)$$
$$y'_x=w^{-1/2}xw_z'.$$We know from convexity of the elliptic curve  that
$w''_{zz}<0$.
$$y''_{xx}=\frac{w'_z}{w^2}-\frac{x^2(w'_z)^2}{w^{3/2}}+\frac{2x^2w''_{zz}}{w^{1/2}}<0 .$$

\bigskip

Assume that $P$ is a convex configuration which is critical for $E$.

$$E=\frac{t}{x} +\frac{1}{y}\ \  \hbox {  implies}$$

$$E'_x=\frac{-t}{x^2}-\frac{y'_x}{y^2}, \ \hbox{ and}$$
$$E''_{xx}=\frac{2t}{x^3}+
\frac{2[y'(x)]^2}{y^3}-
\frac{y''(x)}{y^2}, \ \ \ \hbox{which is positive}.$$

Therefore, each critical point in the convex part is a local
minimum.

So the statement (1) of the proposition is now  straightforward: if
there are two local minima, there should be a local maximum in
between, which is impossible.

(2) For a non-intersecting 4-gon there is a flex strictly increasing
both diagonals.

(3) There should be a maximum point of $E$ in $M(L)$. By what is
proven above, the domain of self-intersecting 4-gons is the only
possibility for it. \qed

\bigskip

\begin{thm}\begin{enumerate}
\item We have a continuous  bijection
$$\varphi:(0,\infty)\rightarrow IntM^C(L)$$ which maps a charge $t$ to the
minimum point of $E$.
    \item The mapping $\varphi$ extends to a bijection
    $\varphi:[0,\infty]\rightarrow M^C(L)$.
\end{enumerate}
\end{thm}
Proof. Point  (1) follows from Proposition \ref{thm4gon}, (1) and
Lemma \ref{UniquetForFourGon}. Point (2)  follows from (1) and the
above remark. \qed

\bigskip

\begin{lemma}Assume that $L$ is a  4-bar linkage, and the charge $t$ is zero.
Then $E$ has exactly two critical points on $M(L)$: one  minimum and
one maximum.
\end{lemma}
Proof. The problem is reduced to finding the critical values of the
diagonal $x$.\qed

\bigskip

The lemma means that if we start at any configuration of any  4-gon
$L$ and put zero charge at the ruling vertex, the gradient flow will
bring us to the global minimum.

\begin{figure}[h]
\centering
\includegraphics[width=12 cm]{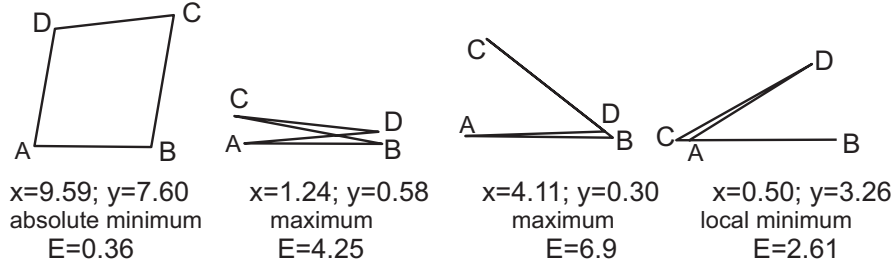}
\caption{Critical configurations from Example
\ref{ExManyCrit}}\label{critical}
\end{figure}

However, for other values of $t$, there might be extra local minima:

\begin{ex}\label{ExManyCrit} For $a=6, \ b=6.5,\ c=6.2, \ d=5.8$ and $t=2$,
numerical computations show that there are  four critical points of
$E$ in total (see Fig. \ref{critical}). The lengths of the diagonals
and the values of $E$ (up to two decimals) are:
\begin{enumerate}
    \item $x=0.50, \ y=3.24,\ E =2.61 $ (local minimum)
    \item $x=4.11, \ y=0.30,\ E =6.90 $ (maximum)
    \item $x=1.24, \ y=0.58,\ E =4.24 $ (local maximum)
    \item $x=9.59, \ y=7.60,\ E =0.36 $ (minimum)
\end{enumerate}
\end{ex}

The above leads us to the following navigating algorithm.
\begin{alg}
Assume we are given a 4-bar linkage $L$, its unknown starting
configuration $P$ and a target convex configuration $P'$. The
following algorithm tells us how to reach $P'$ by altering the
ruling charge $t$.

\begin{enumerate}
    \item Put the ruling charge $t=0$.
    \item Compute the stabilizing charge $t'$ of the configuration $P'$.
    \item Put the ruling charge $t=t'$.
\end{enumerate}
\end{alg}

\bigskip
\begin{rem}

If we allow non-positive charge and aim at non-convex
configurations, the situation becomes more complicated.  Example
\ref{ExManyCrit} shows that the number of (local) minima is greater
than 1.
\end{rem}
\begin{rem}

If we put the same charged 4-bar linkage in 3D, we can skip step 1
in the above algorithm. Indeed,  all critical points of $E$ are
obviously planar configurations. Besides, unlike the 2D case, a
local minimum cannot be self-intersecting: the unfolding flex
increases both of the diagonals. Therefore, in 3D the potential $E$
has exactly one local minimum which is the global minimum.
\end{rem}

\begin{rem}
All above lemmata and Theorem 1 remain valid if we replace the
Coulomb potential  $E$  either  by
\begin{equation} \label{linkcoulombalpha}
E^{\alpha} = \sum \frac{q_iq_j}{(x_{ij})^\alpha}, \ \ \alpha >0,
\end{equation}
or by the limit version of (\ref{linkcoulombalpha})

\begin{equation} \label{linkcoulomblog}
E^{ln} = \sum q_iq_j \ln (x_{ij}).
\end{equation}
\end{rem}

\section{Coulomb control of convex equilateral pentagons}

Let $L$ be a   equilateral  $5$-bar linkage. In our setting we put
two positive charges $s$ and $t$ at a
 pair of non-neighbouring vertices, say, at the third and fifth vertex of
configuration.  The charges of other vertices are constant and equal
to $+1$. Then  the (effective) Coulomb potential $E: M(L)
\rightarrow \mathbb{R}$ is defined by the formula

\begin{equation} \label{coulomb5bar}
E = \frac{1}{x_{14}} + \frac{1}{x_{24}} + \frac{t}{x_{13}} +
\frac{s}{x_{25}} + \frac{st}{x_{35}},
\end{equation}
where $x_{ij}$ are the lengths of diagonals of the configuration
$V$.

We wish to understand whether two non-adjacent charges can provide a
full control (in the same sense as for quadrilateral linkages) on
convex configurations of $L$.
\begin{figure}[h]
\centering
\includegraphics[width=6 cm]{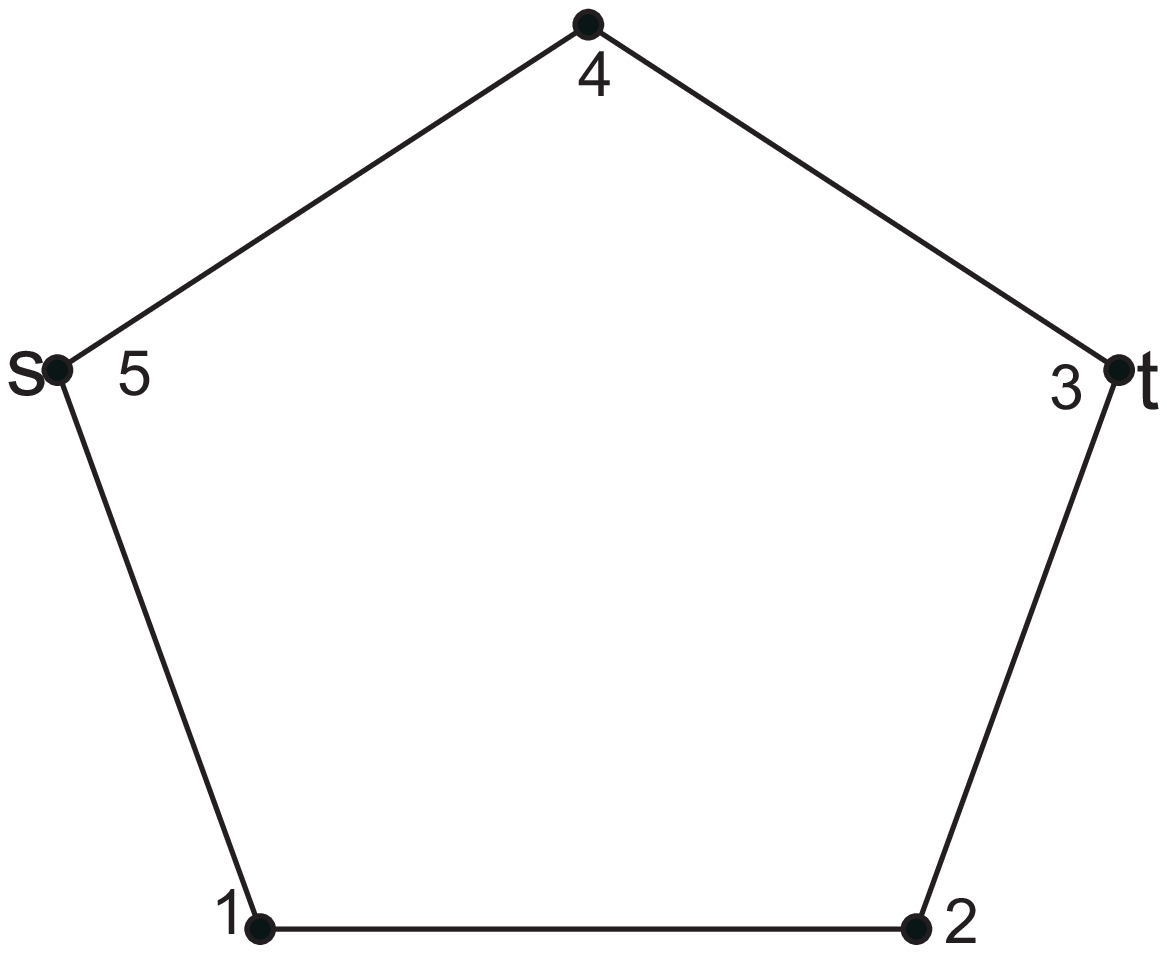}
\caption{}\label{pentagon}
\end{figure}

\begin{lemma}\label{LemmaBoundary} Let $P\in \partial {M^C}$ be a critical point of $E$ for some charges $(s,t)$.
Then $(s,t)$ never belongs to $(0,\infty)\times(0,\infty)$.
\end{lemma}

\begin{figure}[h]
\centering
\includegraphics[width=4 cm]{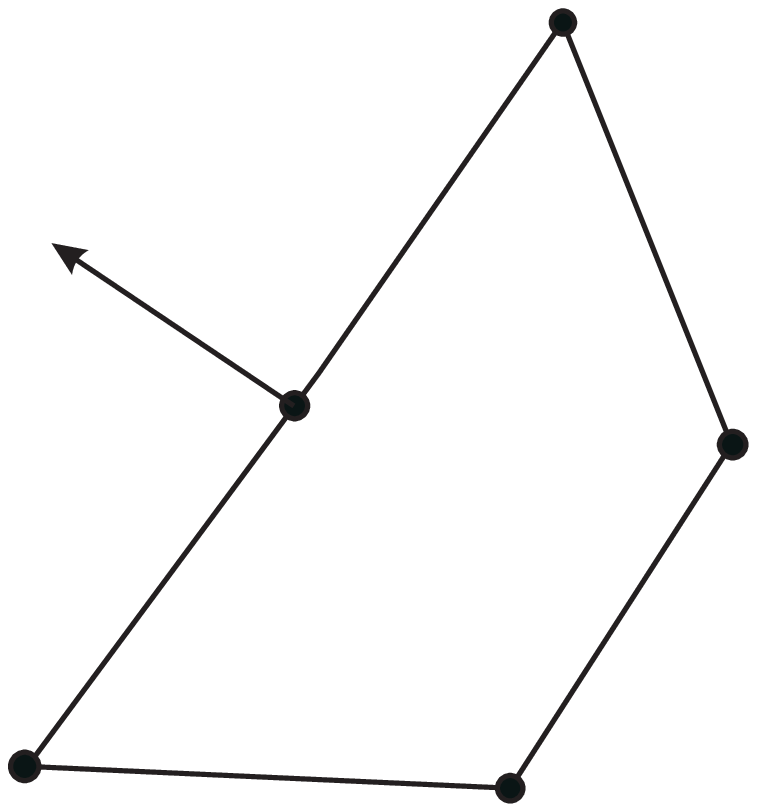}
\caption{}\label{boundary}
\end{figure}

Proof.   Assume the contrary: the charges $s$ and $t$ are both from
$(0,\infty)$. A configuration is a boundary point of $ M^C(L)$
whenever one of the vertices is aligned.
 Then by pulling the aligned vertex  orthogonally to the adjacent edges
 (see Fig. \ref{boundary}) we get an infinitesimal flex which yields a first-order decrease of
$E$.\qed

\begin{figure}[h]
\centering
\includegraphics[width=8 cm]{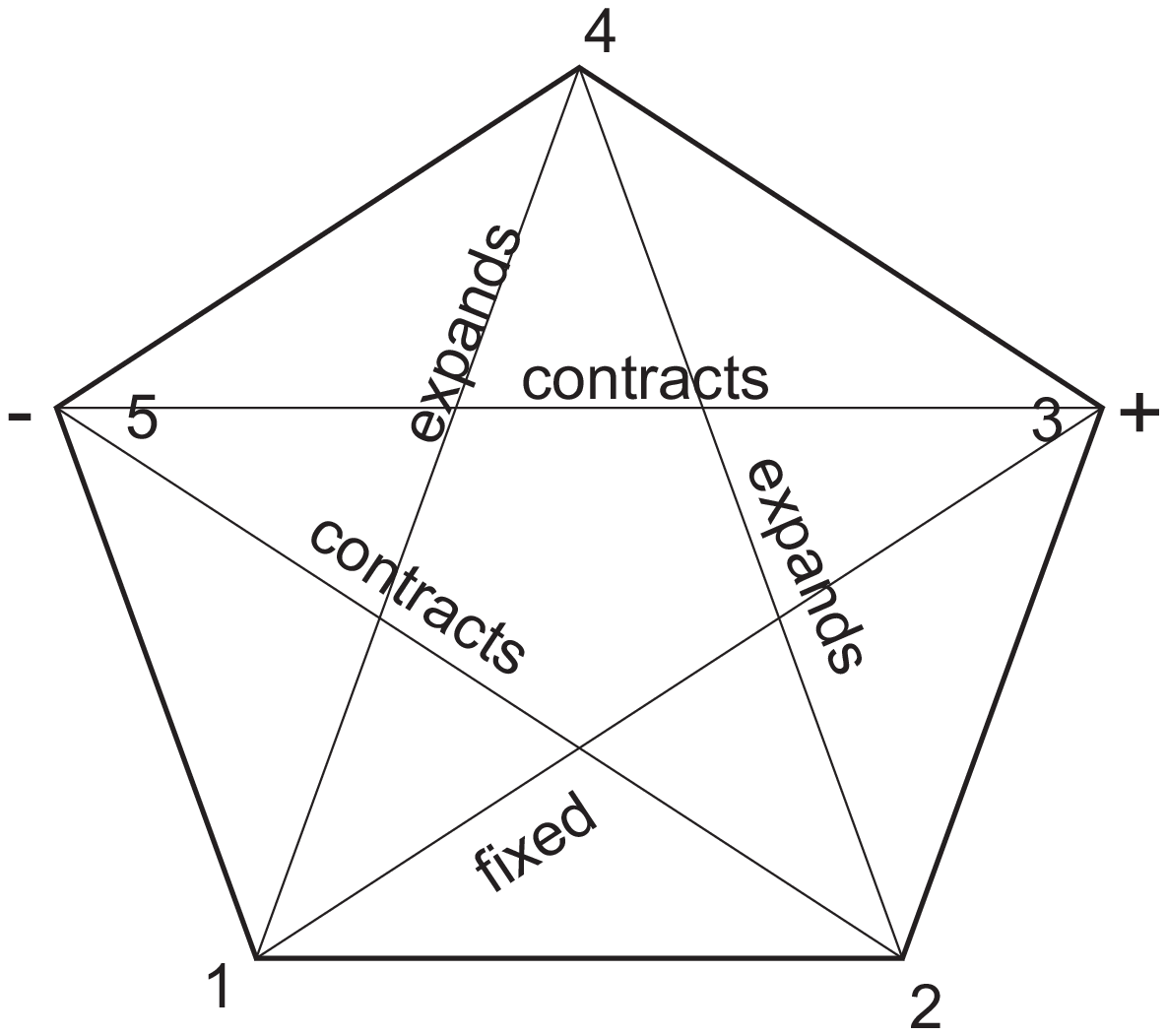}
\caption{}\label{plusminus}
\end{figure}

\begin{lemma}\label{LemmaPlusMinus}
For $s<0,\  t>0$, a convex pentagon is never a critical point of
$E$.
\end{lemma}
Proof. Assume the contrary. Fix $x_{13}$ and compress  $x_{25}$, see
Figure \ref{plusminus}. This  yields an infinitesimal flex such that
$E'<0$.   \qed

\begin{thm}For each strictly convex  configuration $P$ of an  equilateral pentagonal linkage,
there exists exactly one $(s,t) \in (0,\infty )\times (0,\infty )$ such that
 $P$ is a critical point of $E$ for the charges $(s,t)$.

\end{thm}
Proof. We start by reminding that
$$E=1/x_{14}+1/x_{24}+t/x_{13}+s/x_{25}+st/x_{35}.$$

Take the diagonals $x_{35}$ and $x_{13}$ as local coordinates in a
neighborhood of $P$.

The polygon $P$ is a critical point of $E$  means that $dE$
vanishes:
$$-\partial E/ \partial x_{35}=\alpha_1/x_{14}^2+\beta_1/x_{24}^2+ s\gamma_1/x_{25}^2+\frac{st}{x_{35}^2}= 0,
$$
and

$$-\partial E/ \partial x_{13}=\alpha_2/x_{14}^2+\beta_2/x_{24}^2+ t/x_{13}^2+\frac{s\gamma_2}{x_{25}^2 }=
0,
$$

where  $$\alpha_1=\partial x_{14}/
\partial x_{35}, \
 \beta_1=\partial x_{24}/ \partial
x_{35}, \
 \ \gamma_1=\partial x_{25}/ \partial
x_{35}, \hbox{  \  and}$$

$$\alpha_2=\partial x_{14}/
\partial x_{13}, \
 \beta_1=\partial x_{24}/ \partial
x_{13}, \
 \ \gamma_1=\partial x_{25}/ \partial
x_{13}$$

We get a system in two variables $s$ and $t$ which reduces to
the following quadratic equation in $s$:
$$A+Bs+Cs^2=0$$
with
$$A=\frac{\alpha_1}{x_{14}^2}+\frac{\beta_1}{x_{24}^2},$$
$$B=\frac{\gamma_1}{x_{25}^2}-\frac{x_{13}}{x_{35}^2}(\frac{\alpha_2}{x_{14}^2}+\frac{\beta_2}{x_{24}^2}),$$
$$ \hbox{and \ \ \ }C=-\frac{x_{13}^2\gamma_2}{x_{35}^2x_{25}^2}.$$

It is sufficient to prove that $AC$ is negative. Then the equation
has exactly one real positive solution $s$. By Lemma
\ref{LemmaPlusMinus}, $t$ is also positive.

\bigskip
Indeed, $AC$ is negative (and the proof is completed) because of the
following sign analysis.
\begin{enumerate}
    \item Consider the flex of $P$ which fixes  $x_{13}$ and stretches
$x_{35}$.

$x_{14}$ compresses, therefore $\alpha_1=\partial x_{14}/ \partial
x_{35}<0$.

$x_{24}$ compresses, therefore $\beta_1=\partial x_{24}/ \partial
x_{35}<0$.

$x_{25}$ compresses, therefore $\gamma_1=\partial x_{25}/ \partial
x_{35}<0$.

    \item
Similarly,  fixing $x_{35}$ and stretching $x_{13}$, we obtain that
$\alpha_2<0$,
 $\beta_2<0$,
and $\gamma_2<0$.
\end{enumerate}

\begin{figure}[h]
\centering
\includegraphics[width=8 cm]{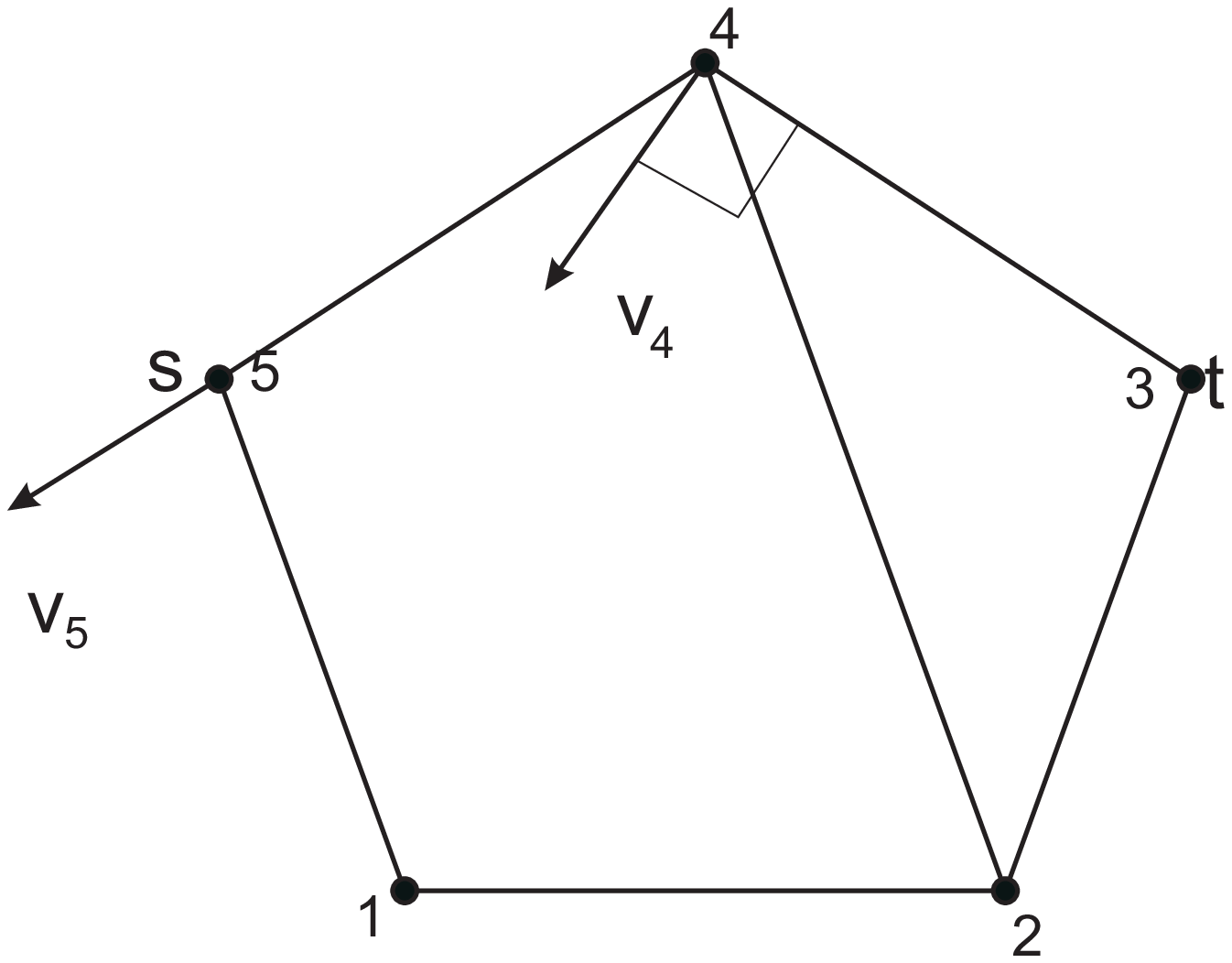}
\caption{}\label{FigVelo}
\end{figure}

Let us give a more detailed explanation why the above inequalities
hold. As an example, let us prove that if the diagonal $x_{13}$ is
fixed, and the diagonal $x_{35}$ gets longer, then $x_{24}$
compresses (see Fig. \ref{FigVelo}).

We can assume that the three points 1,2,3 are fixed and that the
points 4 and 5 are moving with infinitesimal velocities $v_4$ and
$v_5$. The vector $v_4$ is orthogonal to the edge $43$. By
convexity, the angle between $v_4$ and the diagonal $24$ is smaller
than $\pi/2$, which implies the claim. \qed

\begin{rem}
If we try to rule by putting charges at adjacent vertices (say, at
the vertices 1 and 2), the situation gets worse: we reach not all of
the convex polygons. For instance, we will never be able to have
simultaneously the vertices 5 and 3 aligned. By continuity, an
entire neighborhood of such a polygon becomes unreachable.
\end{rem}

\section{Concluding remarks}

Our results suggest several natural problems and conjectures
concerned with the Coulomb control scenario for $n$-bar linkages
with arbitrary $n$. In particular, there is good evidence that, for
equilateral linkages, the complete Coulomb control of convex
configurations is valid if we are permitted to choose controlling
charges at ALL vertices.

Assuming that this conjecture is true, what is the minimal number of
controlling charges sufficient to provide complete Coulomb control?
For dimension reasons, one may await that it should be sufficient to
have $n-3$ controlling charges. Is this amount of controlling
charges indeed sufficient for any $n$-bar linkage?

Several interesting problems in the spirit of the famous Maxwell
conjecture \cite{max} are concerned with estimating the amount and
possible types of critical points of Coulomb potential in the moduli
space of linkage. In particular, what is the exact upper bound for
the number of critical points of Coulomb potential in the planar
moduli space of generic $n$-bar linkage? We conjecture that, for a
generic $4$-bar linkage, the Coulomb potential can have not more
than four critical points and their types are the same as in Example
1.

Analogous problems are also meaningful and interesting for all
potentials of the form (\ref{linkcoulombalpha}),
(\ref{linkcoulomblog}). These and further aspects of the Coulomb
control scenario will be addressed in forthcoming publications of
the authors.


\begin{thebibliography}{99}
\bibitem{conndem}R.Connelly, E.Demaine, \emph{Geometry and topology of plygonal linkages,
Handbook of discrete and computational geometry, 2nd ed.} CRC Press,
Boca Raton, 2004, 197-218.



\bibitem{F}M. Farber, \emph{Invitation to Topological Robotics, Zuerich Lectures in
Advanced Mathematics.} European Mathematical Society (EMS), Zuerich,
2008.

\bibitem{ganosh}
A.Gabrielov, D.Novikov, B.Shapiro, \emph{Mystery of point charges},
Proc. Lond. Math. Soc. 95, 2007, 443-472.


\bibitem{khipan}
G. Khimshiashvili, G. Panina, \emph{Cyclic polygons are critical
points of area}.  Zap. Nauchn. Sem. S.-Peterburg. Otdel. Mat. Inst.
Steklov. (POMI), 2008, 360, 8, 238--245.

\bibitem{khsi}
G.Khimshiashvili, D.Siersma, \emph{Cyclic configurations of planar
multiple penduli}, ICTP Preprint IC/2009/047. 11 p.

\bibitem{kpsz}
G. Khimshiashvili, G. Panina, D. Siersma,  A. Zhukova,
\emph{Critical configurations of planar robot arms,} Centr. Europ.
J. Math., 2013, 11,  3,  519--529.

\bibitem{kundernac} T. Kudernac,N. Ruangsupapichat, M. Parschau, B.
Mac, N. Katsonis, S. Harutyunyan, K.-H. Ernst, B. Feringa, (2011).
\emph{Electrically driven
 directional motion of a four-wheeled molecule on a metal surface}". Nature 479 (7372): 208–11.

\bibitem{max}
J.C.Maxwell, \emph{A Treatise on Electricity and Magnetism,} London,
1853.


\bibitem{GaianeModuliComb}G. Panina, \emph{Moduli space of a planar polygonal linkage:
 a combinatorial description}. arXiv:1209.3241
\end{thebibliography}
\end{document}